\documentclass [12pt,a4paper,reqno]{amsart}
\usepackage{amsmath,amsthm,amscd}
\usepackage{amsfonts,amssymb}
\allowdisplaybreaks

\DeclareMathOperator{\esssup}{ess\,sup}

\newcommand{\be}{\begin{equation}}
\newcommand{\ef}{\end{equation}}
\chardef\bslash=`\\ % p. 424, TeXbook
\newtheorem*{thm*}{Theorem}

\theoremstyle{definition}

\newtheorem*{remark*}{Remarks}
\newtheorem*{defn*}{Definition}

\theoremstyle{remark}
\newcommand{\G}{\Gamma}

\newcommand{\wh}{\widehat}
\newcommand{\bk}{\bigskip}
 \renewcommand{\sectionmark}[1]{}

\newcommand{\qc} {quasiconformal}

\newcommand{\Te} {Teichm\"{u}ller}

\newcommand{\iy} {\infty}
\newcommand{\fc} {\frac}

\usepackage{amsfonts}
\newcommand{\field}[1]{\mathbb{#1}}

\newcommand{\D}{\field D}
\newcommand{\om}{\omega}
\newcommand{\z}{\zeta}
\newcommand{\ov}{\overline}
\newcommand{\vp}{\varphi}
\newcommand{\hC}{\wh{\field{C}}}
\newcommand{\C}{\field{C}}
\newcommand{\R}{\field{R}}
\newcommand{\hR}{\wh{\field{R}}}
\newcommand{\B}{\mathbf{B}}
\newcommand{\T}{\mathbf{T}}

\newcommand{\Bel} {\operatorname{Bel}}

\newcommand{\const}{\operatorname{const}}

\newcommand{\Om} {\Omega}
\newcommand{\vk} {\varkappa}
\newcommand{\x} {\mathbf x}
\renewcommand{\a} {\alpha}

\newcommand{\ld}{\lambda}
\newcommand{\kp}{\kappa}

\begin{document}

\title{A question of K\"{u}hnau}
\author{Samuel L. Krushkal}

\begin{abstract}  
It is well known that the square is not a Strebel's point (i.e., its extremal Beltrami coefficient is not Teichm\"{u}ller). Many years ago, Reiner K\"{u}hnau posed the question:   
does there exist in the case of a``long" rectangle the corresponding holomorphic quadratic differential? 

We prove that the answer is negative for any bounded convex 
quadrilateral and establish a stronger result for rectangles.
\end{abstract}

\date{\today\hskip4mm ({\tt KueQuest.tex})}

\maketitle

\bigskip

{\small {\textbf {2010 Mathematics Subject Classification:} Primary: 30C55, 30C62, 30C70, 30F60}

\medskip

}

\bigskip

\markboth{S. L. Krushkal}{A question of K\"{u}hnau}\pagestyle{headings}

\bigskip
\centerline{\bf 1. EXTREMAL QUASICONFORMALITY}

\bk 
The extremal quasiconformal maps $w(z)$ whose Beltrami coefficients $\mu_w(z) = \overline{\partial} w/\partial w$  
minimize the dilatation $k(w) = \|\mu_w\|_\infty$ play a crucial role in geometric complex analysis and in the Teichm\"{u}ller space theory. 

Consider the disks  $\D = \{z: \ |z| < 1\}, \ \D^* = \{z \in \hC = \C \cup \{\iy\}: \ |z| > 1\}$, and let $\Sigma^0$ be the class of univalent functions $F^\mu(z) = z + b_0 + b_1 z^{-1} + \dots$ on $\D^*$ with quasiconformal extension to $\hC$, so their Beltrami coefficients range over the ball 
$$
\Bel(\D)_1 = \{\mu \in L_\iy(\C): \ \ \mu(z)|\D^* = 0, \ \ 
\|\mu\|_\infty < 1\}. 
$$
Define in this ball the equivalence classes $[\mu]$ and $[F^\mu]$), letting $\mu_1, \ \mu_2$ be equivalent if the corresponding maps $F^{\mu_1}$ and $F^{\mu_2}$ coincide on the unit circle $S^1 = \partial \D$. These classes are in one-to-one correspondence with the Schwarzians $S_{w^\mu}$ on $\D^*$ which fill a bounded domain in the space $\B$ of hyperbolically bounded holomorphic functions on $\D^*$ with norm $\|\vp\| = \sup_{\D^*} (|z|^2 - 1)^2 |\vp(z)|$; it models the {\it universal Teichm\"{u}ller space} $\T$, and the quotient map $
\phi_\T: \ \Bel(\D)_1 \to \T, \ \phi_\T(\mu) = S_{F^\mu}$ 
is a holomorphic split-submersion (has local sections). 

The intrinsic {\it Teichm\"{u}ller metric} on $\T$ is defined by
$$
\tau_\T (\phi_\T (\mu), \phi_\T (\nu)) = \frac{1}{2} \inf \bigl\{\log K \bigl( w^{\mu_*} \circ \bigl(w^{\nu_*} \bigr)^{-1} \bigr) : \ \mu_* \in \phi_\T(\mu), \nu_* \in \phi_\T(\nu) \bigr\},  
$$ 
where $K(w) = (1 + k(w))/(1 - k(w))$. 
It is the integral form of the infinitesimal Finsler metric
$$
F_\T(\phi_\T(\mu), \phi_\T^\prime(\mu) \nu) = \inf
\{\|\nu_*/(1 - |\mu|^2)\|_\iy: \ \phi_\T^\prime(\mu) \nu_* =
\phi_\T^\prime(\mu) \nu\}
$$
on the tangent bundle $\mathcal T \T$ of $\T$.

The Beltrami coefficient $\mu \in \Bel(\D)_1$ is {\it extremal} in its class if
$\|\mu\|_\iy = \inf \{\|\nu\|_\iy: \ \ \phi_\T(\nu) = \phi_\T(\mu)\}$
and {\it infinitesimally extremal} if
$$
\|\mu\|_\iy = \inf \{\|\nu\|_\iy: \ \ \nu \in L_\iy(\D^*), \ \
\phi_\T^\prime(\mathbf 0) \nu = \phi_\T^\prime(\mathbf 0) \mu\}.
$$
Any infinitesimally extremal Beltrami coefficient $\mu$ is globally extremal (and vice versa), and by the basic
Hamilton-Krushkal-Reich-Strebel theorem either of these  extremalities is equivalent to the equality
 \be\label{1}
\|\mu\|_\iy = \inf \{|\langle \mu, \psi\rangle_\D|: \ \ \psi \in A_1(\D): \ \|\psi\| = 1\},
\end{equation}
where $A_1(\D)$ is the subspace of $L_1(\D)$ formed by
holomorphic quadratic differentials $\psi= \psi(z) dz^2$ on $\D$ and the pairing
$$
\langle\mu, \psi\rangle_\D = \iint_\D  \mu(z) \psi(z) dx dy, \quad \mu \in L_\iy(\D), \ \psi \in L_1(\D), \ \ z = x + iy.
$$

Let $F_0 :=F^{\mu_0}$ be an extremal representative of its class $[F_0]$ with dilatation
$$
k(F_0) = \|\mu_0\|_\iy = \sup \{k(F^\mu): F^\mu|S^1 = F_0|S^1\},
$$
and assume that there exists in this class a quasiconformal map $w_1$ whose Beltrami coefficient $\mu_1$ satisfies the
inequality $\esssup_{V_r} |\mu_1(z)| < k(F_0)$ in some annulus $V_r = \{r < |z| < 1\}, \ r > 0$.  Any such $w_1$ is called the {\it frame map} for the class $[w_0]$, and the point of the space $\T$ corresponding to the class $[F_0]$ is called the {\it Strebel point}.

Due to Strebel \cite{St} (see also \cite{EL}), in any  class $[w]$ having a frame map, the extremal map $w_0$ is unique and either conformal or a Teichm\"{u}ller map with Beltrami coefficient $\mu_0 = k |\psi_0|/\psi_0$ on $\D$, defined by an integrable holomorphic quadratic differential $\psi_0$ on $\D$ and a constant $k \in (0, 1)$. 
This is valid, for example, for domains $w(\D)$ bounded by  asymptotically conformal curves (such curves do not have angular points). 

The theorem of Lakic \cite{La}, \cite{GL} yields that the set of Strebel's points is open and dense in $\T$.
Similar results have been established also for arbitrary Riemann surfaces. 

\bigskip
One can define the boundary dilatation $b_x(w)$ of the maps $w$ taking the Beltrami coefficients supported in the neighborhoods of t boundary points $x \in \hR$. The points with $b_x(w) = k(w)$ are called {\it substantial} for $w$ and for its equivalence class (see, e.g., [4, Ch. 17]).

\bigskip\bigskip
\centerline{\bf 2. K\"{U}HNAU'S QUESTION. RESULTS}

\bk
It was established in \cite{Ku2}, \cite{We} that in the case of  the square $P_4^0$ there are different extremal quasiconformal reflections across its boundary (an orientation reversing quasiconformal automorphisms of the sphere $\hC$ preserving 
pointwise the boundary $\partial P_4^0$), all not of Teichm\"{u}ller type. Hence the corresponding class $[\mu]$ defining the maps with $w^\mu(U) = P_4^0$ is not a Strebel point.  

In this connection, Reiner K\"{u}hnau raised many years ago the following 

\bigskip\noindent
{\bf Question}. {\it Does there exist in the case of a "long" 
rectangle the corresponding holomorphic quadratic differential ?} 

\bigskip 
We show that the answer to this question is negative, proving the following stronger results: 
 
\bigskip\noindent
{\bf Theorem 1}. {\it Every bounded convex quadrilateral $P_4$ has a non-Strebel representation in $\T$ (hence, its outer  conformal map $\D^* \to P_4^*$ does not admit an extremal extension of Teichm\"{u}ller type)}. 

\bigskip
For the rectangles, we have more: 

\bigskip\noindent 
{\bf Theorem 2}. {\em There exists a common 
substantial point $z_0 \in S^1$ at which the extremal quasiconformal dilatation $k(F^*)$ of any rectangle is attained.}

\bigskip
The arguments applied in the proof of these theorems extend straightforwardly to arbitrary bounded convex polygons having  
a common outwardly tangent ellipse, in particular to polygons 
obtained by affine transformations of regular polygons. 
So all such polygons have substantial boundary points and thus are not Strebel.

\bk 
We present briefly the needed notions and results underlying the above theorems adapting those to our case; for details see, e.g. \cite{Di}, \cite{EKK}, \cite{GL}, \cite{Kr4}, \cite{Ku2}.  

\bigskip\bigskip
\centerline{\bf 3. A GLIMPSE AT GRUNSKI INEQUALITIES}

\bk 
The proof of the theorem essentially relies on a basic result 
for the Grunski inequalities. We present briefly needed facts from this theory. 

The fundamental Grunsky theorem (extended to multiply connected domains by Milin) states that a holomorphic function $F(z) = z + \const + O(1/z)$ in a neighborhood $U_0$ of the infinite point is extended to a univalent function on the disk $\D^*$ if and only if it satisfies the inequality
$$
\Big\vert \sum\limits_{m,n=1}^\iy \ \sqrt{mn} \ \a_{m n} x_m x_n
\Big\vert \le 1,
$$
where the Grunsky coefficients $\a_{mn}(F)$  are determined by
$$
\log \fc{F(z) - F(\z)}{z - \z} = - \sum\limits_{m,n=1}^\iy \  \a_{mn} z^{-m} \z^{-n}, \quad (z, \z) \in (\D^*)^2,
$$
taking the principal branch of the logarithmic function, and
$\mathbf x = (x_n)$ ranges over the unit sphere $S(l^2)$ of the
Hilbert space $l^2$ of sequences with $\|\mathbf x\|^2 =
\sum\limits_1^\iy |x_n|^2$ (cf. \cite{Gr}). The quantity
$$
\vk(F) = \sup \Big\{\Big\vert \sum\limits_{m,n=1}^\iy \ \sqrt{mn} \
\a_{mn} x_m x_n \Big\vert: \mathbf x = (x_n) \in S(l^2) \Big\} 
$$
is called the {\it Grunsky norm} of the map $F$.

This norm is dominated by the {\it Teichm\"{u}ller norm} $k(F)$ of 
this map, i.e., with the minimal dilatation among quasiconformal extensions of $F$ onto $\D$ (see \cite{Ku1}, \cite{Kr5}); so,
  \be\label{2}
\vk(F) \le k(F) = \tanh \tau_\T(\mathbf 0, S_F),
\end{equation}
where $\tau_\T$ denotes the Teichm\"{u}ller distance on $\T$).  The second norm is intrinsically connected with integrable
holomorphic quadratic differentials on $\D$ (the elements of the subspace $A_1 = A_1(\D)$ of $L_1(\D)$ formed by holomorphic
functions), while the Grunsky norm naturally relates to the {\it abelian} structure determined by the set of quadratic differentials
$$
A_1^2 = \{\psi \in A_1: \ \psi = \om^2\} 
$$
having only zeros of even order on $\D$. It is shown in \cite{Kr1} that such differentials have the form 
 \be\label{3}
\psi(z) = \fc{1}{\pi} \sum\limits_{m+n=2}^{\iy} 
\sqrt{m n} \ x_m x_n z^{m+n-2}
\end{equation}
with $\mathbf x = (x_n) \in l^2$, and $\|\mathbf x\|_{l^2} = \|\psi\|_{A_1}$. 

A crucial point here is that for a generic function $F \in
\Sigma^0$ in (2) the strict inequality $\vk(F) < k(F)$ is valid; moreover, it holds on the (open) dense subset of $\Sigma^0$ in both strong and weak topologies (i.e., in the Teichm\"{u}ller distance and in locally uniform convergence on $D^*$) (see, e.g., \cite{Kr1}, \cite{Kr4}, \cite{Ku2}). 
 
So it is important to know whether for a concrete function 
$F$, we have $\vk(F) = k(F)$. In terms of the pairing
$$
\langle \mu, \psi\rangle_\D = \iint\limits_\D  \mu(z) \psi(z) dx dy,
\quad \mu \in L_\iy(\D), \ \psi \in L_1(\D) \ \ (z = x + iy),
$$ 
such functions are completely characterized by the following 

\bigskip\noindent
{\bf Lemma 1}. \cite{Kr1}, \cite{Kr5} \ {\it For all $F \in \Sigma^0$,
$$
\vk(F) \le k \fc{k + \a(F)}{1 + \a(F) k}, \quad k = k(F),
$$
and $\vk(F) < k$ unless
 \be\label{4}
\a(F) := \sup \ \{|\langle \mu, \psi\rangle_\D|: \ \psi \in A_1^2, \ \|\psi\|_{A_1(\D)} = 1\} = \|\mu\|_\iy;
\end{equation}
the last equality is equivalent to $\vk(F) = k(F)$. Moreover, for small $\|\mu\|_\iy$,
$$
\vk(F) = \sup \ |\langle \mu, \psi\rangle_\D| +
O(\|\mu\|_\iy^2), \quad \|\mu\|_\iy \to 0,
$$
with the same supremum as in (4).

If $\vk(F) = k(F)$ and the equivalence class of $F$ (the collection of maps equal to $F$ on $S^1 = \partial D^*$) is a Strebel point, then the extremal $\mu_0$ in this class is necessarily of the form }
 \be\label{5}
\mu_0 = \|\mu_0\|_\iy |\psi_0|/\psi_0 \ \ \text{with} \ \ \psi_0 \in
A_1^2.
\end{equation}

Geometrically, (4) means the equality of the Carath\'{e}odory and Teichm\"{u}ller distances on the image of the geodesic disk
$\D(\mu_0) = \{t\mu_0 /\|\mu_0\|_\iy: \ t \in \D\}$ in the space
$\T$. For functions $F \in \Sigma^0$ holomorphic in the closed disk $\ov{\D^*}$, the relation (5) was also obtained by a different method in \cite{Ku2}. 

Note also that the Grunsky coefficients $\a_{mn}(F)$ are holomorphic functions of the Schwarzians $\vp = S_F$ on the  universal Teichm\"{u}ller space $\T$, and for each $\x = (x_n) \in S(l^2)$ the series
 \be\label{6}
h_{\mathbf x}(\vp) =
\sum\limits_{m,n=1}^\iy \ \sqrt{mn} \ \a_{mn}(\vp) x_m x_n
\end{equation}
defines a holomorphic map of the space $\T$ into the unit disk $\D$ so that $\vk(F) = \sup_{\x} | h_{\mathbf x}(S_F)|$.

\bigskip\bigskip
\centerline{\bf 4. PROOF OF THEOREM 2}

\bk
We first prove the second theorem answering directly K\"{u}hnau's question and illustrating the main features. 

Let $P_4^0$ be the unit square centered at the origin, and let $\mu_0 \in \Bel(\D)_1$ be an extremal Beltrami coefficient for  $P_4^0$, i.e., for the extremal extension of the outer conformal map $F_{P_0^*}: \ \D^* \to P_0^*$. 

As was mentioned above, this $\mu_0$ is not Teichm\"{u}ller; hence, there exist a point $z_0$ on the unit circle $S^1$ and a degenerating sequence $\{\psi_n\} \subset A_1(\D)$ such that 
$\psi_n(z) \to 0$ locally uniformly on $\D, \ \|\psi_n \|_{A_1(\D)} = 1$, and the boundary dilatation at $z_0$ satisfies  
 \be\label{7}
b_{z_0}(F^{\mu_0}) = k(w^{\mu_0}) = \lim\limits_{p\to\infty} 
|<\mu_0, \psi_p>_\D|.
\end{equation} 
On the other hand, since for the square (as well as for any rectangle) its outer conformal mapping function $F^* = F_{P_0^*}$ has equal Grunsky and Teichm\"{u}ller norms, Lemma 1 and (3) yield that all $\psi_p$ in (7) belong to $A_1^2(\D)$ 
and hence, 
 \be\label{8}
\psi_p(z) = \fc{1}{\pi} \sum\limits_{m+n=2}^{\iy} 
\sqrt{m n} \ x_m^p x_n^p \ z^{m+n-2}, 
\end{equation}  
where $\x^p = (x_1^p, \dots, \ x_n^p, \dots) \in S(l^2)$. 
So from (7) and (8), we have  
 \be\label{9}
b_{z_0}(F^*) = \vk(F^*) = k(F^*) =\lim\limits_{p\to \iy} h_{\x^p}(S_{F^*}).  
\end{equation}
Note also that, by definition of the Grunsky norm, each $h_{\x^p}$ satisfies 
 \be\label{10}
|h_{\x^p}(S_{F^*})| \le \vk(F^*) = k(F^*), 
\quad p = 1, 2, \dots \ .  
\end{equation} 

We shall need also the differential version of the relations (9). Using the variation of $F^\mu(z) = z + b_0 + b_1 z^{-1} + \dots \in \Sigma^0$ with small $\|\mu\|_\iy$ given by 
$$
F^\mu(z) = z -
\fc{1}{\pi} \iint\limits_\D \fc{\mu(w) du dv} {w - z} + O(\|\mu^2\|_\iy), \quad w = u + iv,  
$$
one gets
$$
b_n = \fc{1}{\pi} \iint\limits_\D \mu(w) w^{n-1} du dv +
O(\|\mu^2\|_\iy), \quad n = 1, 2, \dots,
$$
and 
$$ 
\a_{mn}(\mu) = - \pi^{-1}
\iint\limits_\D \mu(z) z^{m+n-2} dx dy + O(\|\mu\|_\iy^2).  
$$ 
Combining this with (3) and (6), one derives that the differential at zero of the corresponding  map $h_\x(t\mu)$ in the direction determined by $\mu$ equals
 \be\label{11}
d h_{\x}(0) \mu = - \fc{1}{\pi} \iint\limits_\D
\mu(z) \sum\limits_{m+n=2}^\iy \sqrt{m n} \ x_m x_n \ z^{m+n-2} dxdy = - \langle \mu, \psi \rangle_\D.  
\end{equation}  

One can assume that the vertices of a rectangle $P_4$ are the points $(\pm 1, 0), (\pm 1, \pm 1 + i a), \ a > 0$, so each $P_4$ is obtained from the square $P_4^0$ by an affine transform  
$$
w = t_1 z + t_2 \ov z
$$ 
with real $t_1, \ t_2$, and use the parameter $t = t_2/t_1 \in (-1, 1)$ measuring the affinity. The map $\mathbf b(t) = S_{F_t^*}: \ (- 1, 1) \to \T$ defines a curve $\G$ in $\T$ whose points represents all rectangles $P_4$. 

Now, using the chain rule for Beltrami coefficients $\mu, \nu$ from the unit ball in $L_\iy(\C)$,  
$$
w^\mu \circ w^\nu = w^{\sigma_\nu (\mu)} \ \ \text{with} \ \ \sigma_\nu (\mu) = (\nu + \mu_1)/(1 + \ov{\nu} \mu_1) 
$$
and $\mu_1(z) = \mu(w^\nu(z)) \ov{w^\nu_z}/w^\nu_z$ (so      
for $\nu$ fixed, $\sigma_\nu (\mu)$ depends holomorphically on $\mu$ in $L_\iy$ norm), we consider the composite maps 
  \be\label{12}
W_t = g^t \circ w^{\mu_0}
\end{equation}  
with complex $t \in \D$. Their Schwarzians $S_{W_t}$ fill in the space $\T$ a holomorphically embedded non-geodesic disk $\Om_0 = \mathbf b(G_0)$ which is the image of a simply connected domain $G_0 \subset \C$ containing the interval $(-1 , 1)$ corresponding to rectangles. 

We now take the restrictions $\wh h_p := h_{\x^p}(S_{F^{\sigma_{\mu_0}(\mu)}})$ of functions $h_{\x^p}$  to the disk $\Om_0$ and apply these restrictions to pull backing the hyperbolic metric $ds = |d \z|/(1 - |\z|^2)$ of $\D$. This provides on the disk $\Om_0$ the conformal metrics $ds = \ld_{\wh h_p}(t) |d t|$ with 
 \be\label{13}
\ld_{\wh h_p}(t) = \wh h_p^* \ld_\D = 
\fc{|\wh h_p^\prime (t)| |d t|}{1 - |\wh h_p(t)|^2} 
\end{equation} 
of Gaussian curvature $- 4$ at noncritical points. 
Consider their upper envelope 
$$
\ld_\iy(t) = \sup \ld_{\wh h_p}(t),
$$
taking the supremum over all $p = 1, 2, \dots$ and all $\mu 
= \mu(t) \in \Bel(\D)_1$ with $\phi_\T(\mu) = S_{F^{\sigma_{\mu_0}(\mu)}} \in \Om_0$, and pick its upper semicontinuous regularization
$$
\wh \ld_\iy(t) =\limsup\limits_{t' \to t} \wh \ld_\iy(t').
$$
The standard arguments (cf. e.g., \cite{Di}, \cite{Kr3}) imply that $\wh \ld_\iy(t)$ is a logarithmically subharmonic metric on $\Om_0$ of the generalized Gaussian curvature at most $- 4$, which means that 
$$
\Delta \log \wh \ld_\iy \ge 4 {\wh \ld_\iy}^2,    
$$ 
where $\Delta$ denotes the {\it generalized} Laplacian  
$$
\Delta u(\z) = 4 \liminf\limits_{r \to 0} \frac{1}{r^2}
\Big\{ \frac{1}{2 \pi} \int_0^{2\pi} u(\z + re^{i \theta}) d
\theta - \ld(\z) \Big\} 
$$ 
(for $u \in C^2$, it coincides with the usual Laplacian). 

Arguing similarly with all functions (6) and taking supremum over all $\x \in S(l^2)$ and $\mu \in \Bel(\D)_1$, one obtains on $\Om_0$ the canonical plurisubharmonic Finsler metric $\ld_\vk$ generated by the Grunsky structure on the space $\T$, also of curvature $\kp(\ld_\vk) \le -4$. 
Both metrics $\ld_\iy$ and $\ld_\vk$ are dominated by the  
infinitesimal Kobayashi metric $\ld_{\mathcal K}$ of the space $\T$ (equal to its Teichm\"{u}ller metric by the Royden-Gardiner  theorem) via  
$$
\ld_\iy \le \ld_\vk \le  \ld_{\mathcal K}.   
$$ 

It is shown in \cite{Kr3} that for any rectangle $P_4$ its outer conformal map $F^*$ has equal Grunsky and Teichm\"{u}ller norms, i.e., $\vk(F^*) = k(F^*)$, and this implies the equality of metrics $\ld_\vk$ and $\ld_{\mathcal K}$ on the whole disk $\Om_0$. 

The above construction yields (see the relations (9)-(11)) that the minimal metric $\ld_\iy$ coincides with 
its dominants $\ld_\vk$ and $\ld_{\mathcal K}$, 
in the base point $t = 0$ representing the square $P_4^0$:   
\be\label{14}
\ld_\iy(0) = \ld_\vk(0) = \ld_{\mathcal K}(0).
\end{equation} 
To compare $\ld_\iy$ with either of two other metrics on the whole disk $\Om_0$, we apply the following basic facts.   
 
\bk\noindent{\bf Lemma 2}. \cite{Kr2} {\it The infinitesimal form $ \ld_{\mathcal K}$ of the Kobayashi-Teichm\"{u}ller  metric on the tangent bundle $\mathcal T(\T)$ of $\T$ is continuous, logarithmically plurisubharmonic in $\vp \in \T$ and has constant holomorphic sectional curvature $\kappa_{\mathcal K}(\vp, v) = - 4$ (hence $\Delta \log \ld_{\mathcal K} = 4 \ld_{\mathcal K}^2$).} 

The global Kobayashi and Teichm\"{u}ller distances are logarithmically plurisubharmonic in each of their variables on $\T \times \T$ (cf. \cite{Kr2}). 

\bigskip\noindent
{\bf Lemma 3}. (Minda's maximum principle \cite{Mi}) 
{\it If a function
$u : \ D \to [- \iy, + \iy)$ is upper semicontinuous in a domain $D \subset \C$ and its generalized Laplacian satisfies the inequality $\Delta u(z) \ge K u(z)$ with some positive constant $K$ at any point $z \in D$, where $u(z) > - \iy$, and if
$$
\limsup\limits_{z \to \z} u(z) \le 0 \ \ \text{for all} \
\z \in \partial D,
$$
then either $u(z) < 0$ for all $z \in D$ or else $u(z) = 0$ for all $z \in \Om$.} 

\bigskip 
Take a sufficiently small neighborhood $U_0$ of the point $t = 0$, and let 
$$
M = \{\sup \ld_{\mathcal K}(t): t \in U_0\};
$$
then in this neighborhood,
$\ld_{\mathcal K}(t) + \ld_\iy(t) \le 2M$ 
and the function
$$
u = \log \fc{\ld_\iy}{\ld_{\mathcal K}} = 
\log \ld_\iy - \log \ld_{\mathcal K}   
$$
satisfies 
$$
\Delta u = 4 \ld_\iy^2 - \ld_{\mathcal K}^2 \ge
8M (\ld_\iy - \ld_{\mathcal K}).
$$
The elementary estimate
$$
M \log(t/s) \ge t - s \quad \text{for} \ \ 0 < s \le t < M
$$
(with equality only for $t = s$) implies that
$$
M \log \fc{\ld_\iy(t)}{\ld_{\mathcal K}(t)} \ge \ld_\iy(t) - \ld_{\mathcal K}(t),
$$
and hence,
$$
\Delta u(t) \ge 4 M^2 u(t).
$$
The equality (14) for the square $P_0$ implies that all  metrics $\ld_\iy, \ \ld_\vk, \ld_{\mathcal K}$ must be equal in the entire disk $\Om_0$ (hence, are  continuous on this disk). 

Now observe that, due to a well-known basic result of  potential theory, a subharmionic function in a domain $\Om \subset \R^n$ can be different (smaller) than its upper semicontinuous  envelope only on a subset $E \subset \Om$ of the capacity zero, hence also of the Lebesgue $n$-measure zero. 

Applying this to our metric (13), one derives from the above construction and  relations (1), (9), (11) that omitting a nowhere dense subset $E \subset \Om_0$, the extremal Beltrami coefficients $\mu_0(\cdot, t) $ corresponding to $t \in \Om_0 \setminus E$ satisfy the equality  
 \begin{equation}\label{15}
\|\mu_0(\cdot, t)\|_\iy = \sup_p |\langle \mu_0(\cdot, t), \psi_p\rangle_\D|.   
\end{equation}
Such an equality can hold only when $\mu_0(\cdot, t)$ are not of Teichm\"{u}ller type (since the Teichm\"{u}ller extremal coefficients do not have the substantial boundary points and cannot attain their dilatation on the degenerated sequences). 

The remaining rectangles $P(t)$ with $t \in E$ also cannot be Strebel points, since by the Lakic theorem the set of such  points is open in the space $\T$. 

It remains to establish the last conclusion of Theorem 2 on existence of a common substantial point $z_0 \in \partial \D$ for all rectangles $P(t), \ t \in \Om_0$. 

Let $\mu_t = \mu_0(\cdot, t)$ and consider the maps $F^{s\mu_t}$ with sufficiently small $|s|$ so that 
$\|S_{ F^{s\mu_t}}\|_\B < 2$. Then 
 \begin{equation}\label{16}
s \langle \mu_t, \psi_p\rangle_\D = 
\langle \nu_{F^{s\mu_t}}, \psi_p \rangle_\D,
\end{equation}
where 
$$
\nu_\vp(z) = \frac{1}{2}(1 - |z|^2)^2 \vp(1/\ov z) 1/{\bar z}^4
$$
means the harmonic Beltrami coefficient of the Ahlfors-Weill extension of the map $F^\mu \in \Sigma^0$ whose Schwarzian  $S_{F^\mu} = \vp$, provided that $\|\vp\|_\B < 2$. These coefficients are connected with the initial extremal coefficients $\mu_t$ via    
$$
\nu_{F^{s\mu_t}} = s \mu_t + \sigma(s, t),
$$
where $\sigma(s, t)$ belong to the annihilating set for $A_1(\D)$, i.e., $\langle \sigma(s, t), \psi \rangle_\D = 0$  for all $\psi \in A_1(\D)$ (this set coincides with $\ker \phi_\T^\prime(\mathbf 0)$). 

Noting that the harmonic coefficients $\nu_\vp$ depend holomorphically on $\vp \in \T$ and the metric $\ld_{\mathcal K}$ is continuous on the space $\T$, one extends by applying (16) the initial equality (15) to $t \in \Om_0 \setminus E$. 
This implies that for each $t \in \Om_0$, the $L_\iy$-length $\|\mu_0(\cdot, t)\|_\iy$ is attained on a subsequence from the initial sequence $\{\psi_p\}$.  This  completes the proof of the theorem.

\bigskip\bigskip
\centerline{\bf 5. PROOF OF THEOREM 1}

\bk 
The idea of the proof of this theorem (without assertion on a common substantial point) is similar to the case of rectangles taking into account that the basic property 
$\vk(F^*) = k(F^*)$ holds for every bounded  rectilinear convex quadrilateral $P_4$.  

We first establish the assertion of Theorem 1 for trapezoids. 
Given a trapezoid $P_4$, one can inscribe an ellipse $\mathcal E$ tangent to all sides of $P_4$. Now take an affine map $g(w)$ moving $\mathcal E$ into a circle. Then in view of affinity, the image $P_4^0 = g(P_4)$ is a rectilinear quadrilateral with a common tangent from the outward circle, and Werner's construction in \cite{We} implies an extremal quasiconformal reflection across the boundary of $P_4^0$ with a variable dilatation function, thus not Teichm\"{u}ller; 
equivalently, $P_4^0$ is not a Strebel point in $\T$. 

The collection of all affine deformations of $P_4^0$, by applying the composed maps (12), generates a holomorphic disk $\Om \subset \T$ representing quadrilaterals. One can repeat for this disk all the above constructions from the proof of Theorem 2 giving the corresponding subharmonic infinitesimal metrics on $\Om$. After a similar comparison, one obtains that the original quadrilateral $P_4$ by itself cannot be obtained from $\D$ by a Teichm\"{u}ller extremal map. 

\bigskip
We now pass to the generic quadrilaterals. It suffices to prove the theorem for bounded quadrilaterals $P_4 = A_1 A_2 A_3 A_4$ (with vertices $A_j$ ordered accordingly  to positive direction of $\partial P_4$) which 
have a line $L$ inside $P_4$ drawn from the vertex $A_1$ parallel to the opposite edge $A_2 A_3$ such that $L$  separates this edge from the remaining vertex $A_4$. 

Fix such a quadrilateral $P_4^0 = A_1^0 A_2^0 A_3^0 A_4^0$ and consider the collection $\mathcal P^0$ of quadrilaterals $P_4 = A_1^0 A_2^0 A_3^0 A_4$ with the same first three vertices and variable $A_4$; the corresponding $A_4$ runs over a subset $E$  of the thrice punctured sphere $\hC \setminus \{A_1^0, A_2^0, A_3^0\}$. One such quadrilateral  is a trapezoid, which we denote by $P_4^* = A_1^0 A_2^0 A_3^0 A_4^*$. We know that its extremal 
Beltrami coefficient $\mu_0$ is not Teichm\"{u}ller and that 
its outer conformal map has equal Grunsky and Teichm\"{u}ller 
norms. 

For any $P_4 \in \mathcal P^0$, the conformal map $F$ of disk $\D^*$ onto the complementary complement $P_4^*$ is represented by the Schwarz-Christoffel integral
 \begin{equation}\label{17}
F(z) = d_1 \int\limits_0^z \prod\limits_1^4(\z - e_j)^{\a_j - 1} \ \frac{d \z}{\z^2} + d_0,
\end{equation}
where $e_j = F^{-1}(A_j) \in S^1, \ \pi \a_j$ is the interior angle at $A_j$ for $P_4^*$, and $d_0, d_1$ are two complex constants. Let $F^0$ denote the conformal map for the complement of $P_4^0$. 

One obtains from the general properties of quasiconformal maps and (17) that the logarithmic derivatives 
$b_F = (\log F^\prime)^\prime = F^{\prime\prime}/F^\prime$
of maps $F$ defining the quadrilaterals $P_4 \in \mathcal P^0$ are (for a fixed $z$)  real analytic functions of $t = A_4$. Passing to their Schwarzians 
$$ 
S_F = b_F^\prime - \fc{1}{2} b_F^2 \in \T,    
$$
one can find a smooth real arc $\G = \mathbf b(t) \subset \T$ containing the point $S_{F^0}$ and the points corresponding to trapezoids; here $\mathbf b$ denotes the map $t = A_4 \to S_F$. 

Since $\T$ is a domain, there is a tubular neighborhood 
containing $\G$; therefore, $\G$ is located on some nonsingular holomorphic disk of the form 
$\Om = \mathbf F(G) \subset \T$, where $G$ again is a simply connected planar domain containing $G$. This disk is not geodesic in the Teichm\"{u}ller-Kobayashi metric on $\T$ and does not pass through the base point of this space, but one can apply to it the same arguments as in the proof of Theorem 2 constructing by (6) the holomorphic maps $h_\x$ for 
the Schwarzians of compositions $W_t = g^t \circ F^*$, where   
the initial map $W_0 = F^*$ is the outer conformal map of a trapezoid $P_4^* \in  \mathcal P^0$ and $g^t$ runs over this collection. 

The restrictions of these maps to the disk $\Om$ determine by pulling back the hyperbolic metric of the unit disk the corresponding conformal metrics of type (13), and one can apply to their upper semicontinuous envelope $\ld_\iy$ the same aruments as in the proof of Theorem 2, getting that 
for any $P_4 \in \mathcal P^0$ its extremal Beltrami coefficients is deternined by a degenerating sequences and hence not Teichm\"{u}ller. This completes the proof of the theorem.

\medskip
{\small\em{ \leftline{Department of Mathematics, Bar-Ilan
University} \leftline{5290002 Ramat-Gan, Israel} \leftline{and
Department of Mathematics, University of Virginia,}
\leftline{Charlottesville, VA 22904-4137, USA}}}


\bigskip
\bigskip
\begin{thebibliography}{EKK}
{\small


\bibitem{Di}
S. Dineen, {\em The Schwarz Lemma},
Clarendon Press, Oxford, 1989. 

\bibitem{EKK}
C.J. Earle, I. Kra  and S.L. Krushkal, {\it Holomorphic motions and Teichm\"{u}ller spaces}, Trans. Amer. Math. Soc. \textbf{944} (1994), 927-948.

\bibitem{EL}
C.J. Earle and Zong Li, {\em Isometrically embedded polydisks in infinite dimensional Teichm\"{u}ller spaces}, J. Geom. Anal. \textbf{9} (1999), 51-71.

\bibitem{GL}
F.P. Gardiner and N. Lakic {\it Quasiconformal Teichm\"{u}ller Theory}, Amer. Math. Soc., Providence, RI, 2000.

\bibitem{Gr}
H. Grunsky, {\it Koeffizientenbediengungen f\"{u}r schlicht
abbildende meromorphe Funktionen}, Math. Z. \textbf{45} (1939),
29-61.

\bibitem{Kr1}
S.L. Krushkal, {\it Grunsky coefficient inequalities,
Carath\'{e}odory metric and extremal \qc \ mappings}, Comment. Math. Helv. \textbf{64} (1989), 650-660.

\bibitem{Kr2}
S.L. Krushkal  {\it Plurisubharmonic features of the \Te \
metric}, Publications de l'Institut Math\'{e}matique-Beograd, Nouvelle s\'{e}rie \textbf{75(89)} (2004), 119-138. 

\bibitem{Kr3}
S.L. Krushkal, {\it Quasireflections, Fredholm eigenvalues and Finsler metrics}, Doklady Mathematics \textbf{69} (2004), 221-224. 

\bibitem{Kr4}
S.L. Krushkal, {\it Strengthened Moser's conjecture, geometry of Grunsky coefficients and Fredholm eigenvalues},
Central European J. Math \textbf{5(3)} (2007), 551-580.

\bibitem{Kr5}
S.L. Krushkal, {\it Strengthened Grunsky and Milin inequalities}, Contemp. Math. \textbf{667} (2016), 159-179. 

\bibitem{Ku1}
R. K\"{u}hnau, {\it Wann sind die Grunskyschen
Koeffizientenbedingungen hinreichend f\"{u}r $Q$-quasikonfor\-me Fortsetzbarkeit}? Comment. Math. Helv. \textbf{61} (1986), 290-307. 

\bibitem{Ku2}
R. K\"{u}hnau, {\it M\"{o}glichst konforme Spiegelung an einer Jordankurve}, Jber. Deutsch. Math. Verein. \textbf{90} (1988),
90-109. 

\bibitem{La}
N. Lakic, {\it Strebel Points}, Lipa's Legacy, Contemporary
Mathematics \textbf{211}, Amer. Math. Soc., Providence, RI, 2001, 417-431.

\bibitem{Mi}
D. Minda, {\it The strong form of Ahlfors' lemma},
Rocky Mountain J. Math. \textbf{17} (1987), 457-461.

\bibitem{St} 
K. Strebel, {\it On the existence of extremal Teichmueller
mappings}, J. Anal. Math. \textbf{30} (1976), 464-480.

\bibitem{We}
S. Werner, {\it Spiegelungskoeffizient und Fredholmscher Eigenwert f\"{u}r gewisse Polygone}, Ann. Acad. Sci. Fenn. Ser. AI. Math., \textbf{22} (1997), 165-186.

}
\end{thebibliography}
\end{document}